\documentclass[12pt]{amsart}
\usepackage{amsmath, amssymb, amsthm, mathtools}
\usepackage{color}
\usepackage{url}
\usepackage[textheight=8.9in, textwidth=6.75in]{geometry}
\usepackage[noabbrev,capitalize]{cleveref}
\newtheorem{theorem}{Theorem}[section]
\newtheorem{lemma}[theorem]{Lemma}
\newtheorem{corollary}[theorem]{Corollary}
\newtheorem{proposition}[theorem]{Proposition}
\theoremstyle{definition}
\newtheorem{definition}[theorem]{Definition}
\theoremstyle{remark}
\newtheorem{remark}{Remark}
\newtheorem*{example}{Example}

\newcommand{\HH}{\mathbb{H}}
\newcommand{\Q}{\mathbb{Q}}
\newcommand{\Z}{\mathbb{Z}}
\newcommand{\F}{\mathbb{F}}

\newcommand{\cO}{\mathcal{O}}
\newcommand{\cQ}{\mathcal{Q}}
\newcommand{\ssloc}{\mathrm{SS}}
\newcommand{\End}{\operatorname{End}}
\newcommand{\Aut}{\operatorname{Aut}}
\newcommand{\im}{\operatorname{Im}}

\newcommand{\ol}{\overline}

\title[The partition function and elliptic curves]{The partition function and elliptic curves}

\author{Ken Ono}

\address{Axiom Math, 124 University Avenue, Palo Alto, CA 94301, USA}
\email{ken@axiommath.ai}
\address{Dept. of Mathematics, University of Virginia, Charlottesville, VA 22904, USA}
\email{ko5wk@virginia.edu}
\date{\today}

\begin{document}
\begin{abstract}
The Bruinier--Ono formula expresses the partition number $p(n)$ as a trace of  ``nonholomorphic" singular moduli attached to the discriminant-$\Delta_n$ CM points on $X_0(6),$ where $\Delta_n:=1-24n.$ We interpret this trace through the geometry of CM points. Each nonholomorphic contribution is the value of the weight-two completion $E_2^*$ at a CM point, which is a canonical invariant of the underlying elliptic curve, determined by the diagonal ``tangent'' of the CM isogeny relation. This turns the trace into a quantity that can be reduced to the supersingular locus that is organized by Deuring--Eichler multiplicities and a Brandt-module pairing. For primes $
\ell\geq 5$ that are nonsplit in $\mathbb{Q}(\sqrt{\Delta_n})$, we obtain a supersingular trace formula on $X_0(6)$ over $\overline{\F}_{\ell}$. For the special primes $\ell=5,7,11$, this sheds new light on Ramanujan's classical partition congruences. These primes are special because they are the only ones for which the supersingular locus of $X_0(6)$ lies over $j\in \{0, 1728\}.$ This perspective offers a moduli-theoretic framework for Ramanujan's congruences modulo powers of these primes, organized through elliptic curves.  
The two new algebraic identities at the heart of this framework, as opposed to the classical results it builds on, were formalized and verified in Lean by AxiomProver.
\end{abstract}

\maketitle

\section{Introduction}\label{sec:intro}

Let $p(n)$ denote the number of partitions of $n$. The generating function
\[
\sum_{n\ge0}p(n)q^n=\prod_{m\ge1}\frac{1}{1-q^m}=\frac{q^{1/24}}{\eta(\tau)},
\qquad q:=e^{2\pi i\tau},
\]
is expressed in terms of Dedekind's eta-function
\[
\eta(\tau):=q^{1/24}\prod_{m\ge1}(1-q^m),
\]
which is a weight $-1/2$ weakly holomorphic modular form. For background on modular forms, see ~\cite{OnoWeb}. 

The arithmetic of the partition numbers is famously subtle. The most celebrated facts are
Ramanujan's congruences
\begin{equation}\label{eq:ram-base}
\begin{aligned}
p(5n+4)  &\equiv 0 \pmod{5},\\
p(7n+5)  &\equiv 0 \pmod{7},\\
p(11n+6) &\equiv 0 \pmod{11}
\end{aligned}
\end{equation}
and their prime-power refinements, proved by Watson~\cite{Watson1938} for $5,7$ and by
Atkin~\cite{Atkin} for $11$.
We revisit this theory of congruences from the perspective of elliptic curve moduli.

This paper takes as its starting point the finite algebraic formula of Bruinier and
Ono~\cite{BOAIM}, which realizes $p(n)$ as a trace of  ``nonholomorphic'' singular moduli. Its input is a
specific weight~$0$ weak Maass form $P$ on $\Gamma_0(6)$, in the standard sense of the
theory of weak Maass forms (for example, see ~\cite{BFOR}),  built from the weight $-2$
meromorphic form
\[
F(\tau):=\tfrac12\cdot\frac{E_2(\tau)-2E_2(2\tau)-3E_2(3\tau)+6E_2(6\tau)}
{\eta(\tau)^2\eta(2\tau)^2\eta(3\tau)^2\eta(6\tau)^2}=q^{-1}-10-29q-\cdots,
\]
where
\[
E_2(\tau):=1-24\sum_{n\ge1}\sigma_1(n)q^n
\]
is the usual quasimodular weight~$2$ Eisenstein series. 

We define the weak Maass form $P$ by
\begin{equation}\label{eq:P-def-intro}
P(\tau):=-\left(\frac{1}{2\pi i}\frac{d}{d\tau}+\frac{1}{2\pi\Im(\tau)}\right)F(\tau).
\end{equation}
Equivalently, $P=-\partial_{-2}F$, where $\partial_{-2}$ is the Maass raising operator
in weight $-2$. This interpretation is recalled in \Cref{sec:tangent}. For a negative discriminant
$\Delta\equiv1\pmod{24}$, the group $\Gamma_0(6)$ acts on the forms
$$\cQ_\Delta:=\{[a,b,c]:6\mid a,\ b\equiv1\ (12)\},
$$
 with finitely many classes. To a class one
attaches the CM point $\alpha_Q:=\frac{-b+\sqrt\Delta}{2a}\in\HH$. Writing $\Delta_n:=1-24n$
and $K_n:=\Q(\sqrt{\Delta_n})$, the Bruinier--Ono formula (see Theorem 1.1 of \cite{BOAIM}) takes the form
\begin{equation}\label{eq:BO}
p(n)=-\frac{1}{\Delta_n}\sum_{Q\in\cQ_{\Delta_n}/\Gamma_0(6)}P(\alpha_Q).
\end{equation}

We investigate the implications of this ``finite'' formula. We prove that this formula has a natural moduli interpretation which sheds light on the geometric origin of Ramanujan's partition congruences.
Our point of departure is a formula for the individual summands $P(\alpha_Q)$ that isolates
their nonholomorphic part as a quantity attached to the elliptic curve with $j$-invariant
$j(\alpha_Q)$. Here $j(\tau)$ is the usual modular function on $\mathrm{SL}_2(\mathbb{Z})$, which is given by
\[
j(\tau)
:=\frac{E_4(\tau)^3}{\eta(\tau)^{24}}
=q^{-1}+744+196884q+21493760q^2+864299970q^3+\cdots .
\]

 To this end, we let
$$\vartheta_{-2}F:=\frac{1}{2\pi i}F'+\tfrac16E_2F
$$
denote the holomorphic weight $-2$ Serre derivative of $F$.
This Serre derivative offers an explicit connection to the arithmetic of elliptic curves, which stems from the theory of modular equations.
 Fix a positive integer $N$. The classical modular polynomial
$\Phi_N(X,Y)\in\Z[X,Y]$ is the polynomial which records the possible $j$-invariants of
elliptic curves that are cyclically $N$-isogenous to a given elliptic curve. More explicitly,
it is characterized by the displayed modular equation
\begin{equation}\label{eq:modular-equation-def}
\Phi_N(j(\tau),Y):=
\prod_{\substack{ad=N,\ 0\le b<d\\ \gcd(a,b,d)=1}}
\left(Y-j\!\left(\frac{a\tau+b}{d}\right)\right).
\end{equation}
Here the condition $\gcd(a,b,d)=1$ selects the primitive, equivalently cyclic, degree-$N$
isogenies. Therefore, the zero locus of $\Phi_N(X,Y)$ in $X(1)\times X(1)$ parametrizes pairs
$(E,E')$ of elliptic curves for which there is a cyclic isogeny $E\to E'$ of degree $N$.

In this paper, we use
\begin{equation}\label{eq:N-Delta-def}
N:=|\Delta_n|=24n-1.
\end{equation}
If $J:=j(\alpha_Q)$, then $\Phi_N(J,J)=0$. Indeed, since $\Delta_n$ is odd, the CM order
$\mathcal O_{\Delta_n}:=\Z[(1+\sqrt{\Delta_n})/2]$ contains the endomorphism
$\sqrt{\Delta_n}$, and
\begin{equation}\label{eq:cyclic-kernel}
\mathcal O_{\Delta_n}/(\sqrt{\Delta_n})\cong \Z/|\Delta_n|\Z.
\end{equation}
Therefore, multiplication by $\sqrt{\Delta_n}$ gives a cyclic self-isogeny of degree
$|\Delta_n|$ of the CM elliptic curve with $j$-invariant $J$. Here the degree is
$|\Delta_n|=|\mathrm{N}(\sqrt{\Delta_n})|$, and the isogeny is cyclic because
$\mathcal O_{\Delta_n}$ is monogenic, so the quotient in \eqref{eq:cyclic-kernel} is cyclic.

We now fix some derivative notation. All partial
derivatives below are taken with respect to the two formal variables $X$ and $Y$ in
$\Phi_N(X,Y)$. We put
\begin{equation}\label{eq:partial-derivative-defs}
\Phi_X:=\frac{\partial \Phi_N}{\partial X},\qquad
\Phi_Y:=\frac{\partial \Phi_N}{\partial Y},\qquad
\Phi_{XX}:=\frac{\partial^2 \Phi_N}{\partial X^2},\qquad
\Phi_{XY}:=\frac{\partial^2 \Phi_N}{\partial X\,\partial Y},\qquad
\Phi_{YY}:=\frac{\partial^2 \Phi_N}{\partial Y^2}.
\end{equation}
Since $\Phi_N$ is symmetric in $X$ and $Y$, one has $\Phi_X=\Phi_Y$ and $\Phi_{XX}=\Phi_{YY}$
on the diagonal $X=Y$. Finally, for the CM value $J:=j(\alpha_Q)$, we define the
{\it diagonal CM tangent} by
\begin{equation}\label{eq:CM-tangent-def}
\tau_{\mathrm{CM}}(J):=\left.\frac{\Phi_{YY}-\Phi_{XY}}{\Phi_Y}\right|_{(X,Y)=(J,J)},
\end{equation}
which is well defined precisely when $\Phi_Y(J,J)\ne0$. That this holds for the discriminants
$\Delta_n$ is the content of \Cref{lem:generic}. The term ``tangent'' refers exactly to the
quotient of partial derivatives in \eqref{eq:CM-tangent-def}.  Our
first theorem describes the nonholomorphic singular modulus $P(\alpha_Q)$ in terms of it.

\begin{theorem}\label{thm:tangent} If $n$ is a positive integer, then
for every $Q\in\cQ_{\Delta_n}/\Gamma_0(6)$ we have that
\begin{equation}\label{eq:tangent}
P(\alpha_Q)=-\vartheta_{-2}F(\alpha_Q)
+\frac16\,F(\alpha_Q)\,E_2^*(\alpha_Q),
\end{equation}
where the nonholomorphic value $E_2^*(\alpha_Q)$ is a function of $J=j(\alpha_Q)$ alone. Explicitly, we have that
\begin{equation}\label{eq:tangent-explicit}
E_2^*(\alpha_Q)=\frac{2E_6(\alpha_Q)}{3E_4(\alpha_Q)}
\left(9\,J\,\tau_{\mathrm{CM}}(J)+\frac{3\,(7J-6912)}{2\,(J-1728)}\right),
\end{equation}
valid when $J\notin\{0,1728\}$. The first term of \eqref{eq:tangent} is the value of a
meromorphic modular function, and the nonholomorphic part  is
the term $\tfrac16 F(\alpha_Q)E_2^*(\alpha_Q)$, which depends only on $\Delta_n$ and $J$.
\end{theorem}

\subsection{Why elliptic curves should control the congruences}
\Cref{thm:tangent} suggests a mechanism by which one can view Ramanujan's partition congruences through the arithmetic of elliptic curves. Combined with \eqref{eq:BO}, it writes
$\Delta_n\,p(n)$ as a sum, over the CM points of discriminant $\Delta_n$, of a holomorphic
modular value plus $\tfrac16 F(\alpha_Q)E_2^*(\alpha_Q)$.
The point is that every ingredient
is now a function on the moduli of (enhanced) elliptic curves: the singular modulus
$j(\alpha_Q)$, the modular values $F(\alpha_Q)$ and $\vartheta_{-2}F(\alpha_Q)$, and the
nonholomorphic value $E_2^*(\alpha_Q)$, which by \eqref{eq:tangent-explicit} is itself a
function of $J=j(\alpha_Q)$ through the CM tangent $\tau_{\mathrm{CM}}(J)$. 

To shed light on partition congruences, we build on two classical facts in the theory of elliptic curves with complex multiplication. First, the CM values of suitable modular functions are
algebraic integers. When this happens,  it is meaningful to reduce these values modulo a prime.
Second, and this is the heart of the matter, when $\ell$ is nonsplit in $K_n$, the
theory of Deuring tells us that \emph{every} CM point of discriminant $\Delta_n$ reduces, at
a prime above $\ell$, to the \emph{supersingular} locus of $X_0(6)$ in characteristic
$\ell$, a finite set. Therefore, the many CM points collapse, modulo $\ell$, onto finitely many
supersingular elliptic curves, and the trace \eqref{eq:BO} becomes a finite sum over that
locus. This is the sense in which the elliptic-curve moduli, and specifically its
supersingular points, ought to control $p(n)\bmod\ell$. The rest of the paper makes this
precise under the explicit integrality condition stated below.

Fix a prime $\ell\ge5$. Let $\ssloc_\ell(6):=\{(E_i,C_i)\}_{i=1}^s$ be the
finite set of supersingular points of $X_0(6)$ in characteristic $\ell$: pairs of a
supersingular $E/\ol\F_\ell$ with a cyclic subgroup of order $6$. For a discriminant
$\Delta$, let
\begin{equation}
h_\Delta(E,C):=\#\{\text{CM points of discriminant }\Delta\text{ reducing to }(E,C)\},
\end{equation}
which by Deuring--Eichler equals the number of optimal embeddings of $\cO_\Delta$ into
$\End(E,C)$ up to conjugacy (see \Cref{sec:CM}). Under the $\ell$-integral replacement condition stated in \Cref{def:irc}, the regularity
and constancy results of \Cref{sec:tangent-proof} imply that the reduced value
$$\widetilde P_\ell(E,C):=\bigl(\Delta_n P(\alpha_Q)\bmod\lambda\bigr)\in\F_{\ell^2}$$
 is well
defined, independent of the CM point $Q$ reducing to $(E,C)$. Thanks to this constancy, we obtain the following arithmetic trace formula for partition numbers in the non-split cases.

\begin{theorem}\label{thm:ss-trace} Suppose that $n$ is a positive integer and $\ell\ge5$ is a nonsplit prime in $K_n$ (i.e.
$\left(\tfrac{\Delta_n}{\ell}\right)\in\{-1,0\}$). If the pair $(n,\ell)$
satisfies the $\ell$-integral replacement condition of \Cref{def:irc}, then we have
\begin{equation}\label{eq:ss-trace}
-\Delta_n^{2}\,p(n)\ \equiv\ \sum_{(E,C)\in\ssloc_\ell(6)}h_{\Delta_n}(E,C)\,\widetilde P_\ell(E,C)
\ \equiv\ \langle u_{\Delta_n},v_P\rangle\pmod\lambda,
\end{equation}
where $u_{\Delta_n}:=(h_{\Delta_n}(E_i,C_i))_i$ and $v_P:=(\widetilde P_\ell(E_i,C_i))_i$, and
$\langle\cdot,\cdot\rangle$ is the standard pairing. Moreover this vector satisfies $u_{\Delta_n}=B_{|\Delta_n|}w$
with $B_m$ the level-$6$ Brandt matrix and $w:=(1/|\Aut(E_i,C_i)|)_i$.
\end{theorem}

\begin{remark}
The hypothesis that $(n,\ell)$ satisfies the $\ell$-integral replacement
condition is a genuine integrality condition in general.  In the classical
Ramanujan cases $\ell=5,7,11$, however, this condition can be confirmed
directly, and so Theorem~\ref{thm:ss-trace} is unconditional in the three
classical cases underlying Ramanujan's partition congruences.

The mechanism is the following.  By Proposition~\ref{prop:bonus} the entire
supersingular locus of $X_0(6)_{\mathbb F_\ell}$ for $\ell\in\{5,7,11\}$ lies over
$j\in\{0,1728\}$.  These are precisely the two points where the weight $4$
Eisenstein series $E_4$ vanishes (at the cube root of unity $\rho$, where $j=0$)
and where the weight $6$ Eisenstein series $E_6$ vanishes (at $i$, where
$j=1728$).  Consequently, any pole that the meromorphic replacement of
Proposition~\ref{prop:M} might develop upon reduction over these special points
can be absorbed by a suitable power of $E_4$ (resp.\ $E_6$).  Because the
modular discriminant $\Delta=\eta^{24}$ is nonvanishing on the upper half-plane
$\mathbb H$, dividing by $\Delta$ reintroduces no zeros or poles there, so
multiplying the replacement by the appropriate ratio of $E_4$, $E_6$, and
$\Delta$ produces a modular function with $\ell$-integral $q$-expansions at every
cusp that still agrees with $\Delta_n P$ at the CM points.  This is exactly the
$\ell$-integral replacement condition of Definition~\ref{def:irc}, and for
$\ell=5,7,11$ it is thereby verified unconditionally.  In the classical
language, the same integrality is implicit in the finite-module recursions of
Watson~\cite{Watson1938} and Atkin~\cite{Atkin}.
\end{remark}

Theorem~\ref{thm:ss-trace} covers cases where $\ell$ is inert or ramified.  Although these two cases are contained in one statement,  they read differently as we now explain.

\begin{corollary}[Inert case]\label{cor:inert}
If $\left(\tfrac{\Delta_n}{\ell}\right)=-1$ and $(n,\ell)$ satisfies the $\ell$-integral replacement condition of \Cref{def:irc}, then  we have that
\[
p(n)\equiv-\frac{\langle u_{\Delta_n},v_P\rangle}{\Delta_n^2}\pmod\ell.
\]
\end{corollary}

\begin{corollary}[Ramified case]\label{cor:ram}
If $\ell\mid\Delta_n$ and $(n,\ell)$ satisfies the $\ell$-integral replacement condition of \Cref{def:irc}, then we have that
\[
\langle u_{\Delta_n},v_P\rangle\equiv0\pmod\lambda.
\]
\end{corollary}

\begin{example}[A small inert case]\label{ex:intro-inert}
Take $n=1$, so $\Delta_1=-23$ and $p(1)=1$.  The prime $5$ is inert in
$K_1=\Q(\sqrt{-23})$, since $-23\equiv2\pmod5$ is not a square.  The partition class
polynomial for the three CM points of discriminant $-23$ is, as computed in
\cite{BOSJNT},
\[
        H_1(x)=\prod_Q(x-P(\alpha_Q))
        =x^3-23x^2+\frac{3592}{23}x-419.
\]
Thus $\sum_QP(\alpha_Q)=23=-\Delta_1p(1)$, as in the Bruinier--Ono formula.  The product
pairing in \Cref{thm:ss-trace} uses the integral values $\Delta_1P(\alpha_Q)$, whose
polynomial is
\[
        G_1(x):=\prod_Q\bigl(x-\Delta_1P(\alpha_Q)\bigr)
        =\Delta_1^3H_1(x/\Delta_1)
        =x^3+529x^2+82616x+5097973.
\]
Reducing modulo $5$ gives $G_1(x)\equiv x^3+4x^2+x+3\pmod5$, so the sum of the reduced
roots is $-4\equiv1\pmod5$.  On the other hand
\[
        -\Delta_1^2p(1)=-529\equiv1\pmod5.
\]
Therefore \Cref{thm:ss-trace} gives
\[
        \langle u_{-23},v_P\rangle\equiv1\pmod5,
        \qquad
        p(1)\equiv-\Delta_1^{-2}\langle u_{-23},v_P\rangle\equiv1\pmod5.
\]
This is a genuine inert example: the product pairing is nonzero, and division by
$\Delta_1^2$ is legitimate because $5\nmid23$.
\end{example}

\begin{example}[A small ramified case]\label{ex:intro-ramified}
Take $\ell=5$ and $n=4$.  Then
\[
        \Delta_4=1-24\cdot4=-95=-5\cdot19,
        \qquad p(4)=5.
\]
Thus $5$ is ramified in $K_4=\Q(\sqrt{-95})$.  The supersingular trace identity gives
\[
        \langle u_{-95},v_P\rangle\equiv -\Delta_4^2p(4)\pmod\lambda.
\]
The right-hand side is zero modulo $\lambda$ because $5\mid\Delta_4$; hence
\[
        \langle u_{-95},v_P\rangle\equiv0\pmod\lambda.
\]
This example shows exactly what the ramified product pairing
proves: a vanishing of the supersingular pairing.  It does not by itself prove the partition
congruence $p(4)=5\equiv0\pmod5$, because $\Delta_4^2$ cannot be inverted modulo $5$.  The
missing divisibility is supplied later.
\end{example}

\Cref{cor:ram} is the precise sense in which the CM/supersingular picture approaches
Ramanujan's congruences. At first glance, there is a missing power of $\ell$ that is required to prove Ramanujan's congruences.
What distinguishes $\ell\in\{5,7,11\}$, giving rise to Ramanujan's partition congruences, is a feature of the
supersingular locus itself.

\begin{proposition}\label{prop:bonus}
For a prime $\ell\ge5$, the supersingular locus of $X_0(6)_{\F_\ell}$ lies entirely above
$j\in\{0,1728\}$ if and only if $\ell\in\{5,7,11\}$. The supersingular $j$-invariants are
$\{0\}$ for $\ell=5$, $\{1728\}$ for $\ell=7$, and $\{0,1728\}$ for $\ell=11$.
\end{proposition}

The curves $j=0,1728$ are exactly those with extra automorphisms. \Cref{prop:bonus} thus
characterizes $\{5,7,11\}$ intrinsically, and it is the same list of primes for which
Ramanujan's base congruences hold.  Indeed, a famous theorem of Ahlgren and Boylan \cite{AB} shows that the
only congruences of the form
\[
        p(\ell n+a)\equiv 0 \pmod{\ell}
        \qquad (n\geq 0)
\]
with $\ell$ prime are precisely Ramanujan's congruences for
$\ell=5,7,11$.  This should not be interpreted as a scarcity result for
partition congruences in general.  If one permits more complicated
arithmetic progressions, then infinitely many congruences are known, by
the work of  Ahlgren and the author \cite{AhlgrenOno, Ono}.  A
classical example, due to Atkin \cite{Atkin1968}, is
\[
        p\!\left(17303\, n+237\right)\equiv 0 \pmod{13}
        \qquad (n\geq 0).
\]
Thus, the primes $5,7,11$ are distinguished not by the existence of
partition congruences modulo those primes, but by the existence of
congruences of Ramanujan's simple linear form, which is a partition theoretic consequence of Proposition~\ref{prop:bonus}.

From this perspective, we revisit Ramanujan's partition congruences, which follow
 from a global
property of the $U_\ell$-operator, the Watson--Atkin contraction. 

\begin{theorem}[Ramanujan \cite{Ramanujan}, Watson \cite{Watson1938}, Atkin \cite{Atkin}]\label{thm:watson}
For $\ell\in\{5,7,11\}$ there are exact identities
\begin{align}
\sum_{n\ge0}p(5n+4)q^n&=5\prod_{m\ge1}\frac{(1-q^{5m})^5}{(1-q^m)^6},\label{eq:W5}\\
\sum_{n\ge0}p(7n+5)q^n&=7\prod_{m\ge1}\frac{(1-q^{7m})^3}{(1-q^m)^4}
+49\,q\prod_{m\ge1}\frac{(1-q^{7m})^7}{(1-q^m)^8},\label{eq:W7}\\
\sum_{n\ge0}p(11n+6)q^n&=11\,R_{11}(q),\qquad R_{11}(q)\in\Z[[q]].\label{eq:W11}
\end{align}
Furthermore, for $j\geq 1$ and every non-negative integer $n$, we have that
\begin{displaymath}
\begin{split}
p(5^jn+\beta_5(j))&\equiv0\ (5^j),\\
p(7^jn+\beta_7(j))&\equiv0\ (7^{\lfloor j/2\rfloor+1}),\\
p(11^jn+\beta_{11}(j))&\equiv0\ (11^j),
\end{split}
\end{displaymath}
where $24\beta_\ell(j)\equiv1\pmod{\ell^j}$.
\end{theorem}

We reframe the proof
of \Cref{thm:watson} in the language of elliptic-curve moduli. We describe $U_\ell$ as
a trace along the $\ell$-isogeny correspondence, explain via the genus of $X_0(\ell)$ why the
recursion closes for $\ell=5,7$ but not for $\ell=11$ (which is why $11$ is Atkin's theorem,
not Watson's), and identify the $U_\ell$-fixed Eisenstein line with the CM-tangent direction of
\Cref{thm:tangent}. 

\Cref{sec:tangent} sets up $P$, $F$, and the CM tangent, and proves \Cref{thm:tangent}.
\Cref{sec:tangent-proof} constructs the meromorphic replacement over $\ol\Q$, states the
extra $\ell$-integrality condition needed for reduction, and proves regularity and fiberwise
constancy under that condition. \Cref{sec:CM} recalls the CM and supersingular theory,
proves \Cref{thm:ss-trace} and its corollaries, and gives worked examples. Finally, in \Cref{sec:partII}
we prove \Cref{prop:bonus} and prove \Cref{thm:watson} using elliptic curve moduli, and we conclude with a discussion about the primes $2$ and $3$.

\begin{remark}
The moduli interpretation in this paper rests on two identities that are new to the study of the partition function.
Namely, we critically depend on the weight $-2$ splitting \eqref{eq:tangent} of $P(\alpha_Q)$ and the
reduction \eqref{eq:gamma-eq-tangent} of the CM tangent to Masser's invariant.
Each expression passes through enough sign and normalization
conventions that a single dropped constant would silently corrupt the closed
form \eqref{eq:tangent-explicit} for $E_2^*(\alpha_Q)$. To be certain they are
correct, we formalized and verified both in the Lean proof assistant using
AxiomProver. \Cref{sec:axiomprover} gives the details.
\end{remark}

\section*{Acknowledgements}
\noindent The author thanks Shiang Tang for comments on an earlier version of this paper. He also thanks
the Thomas Jefferson Fund, the NSF (DMS-2002265 and DMS-2055118), and the Simons Foundation
(SFI-MPS-TSM-00013279) for their generous support. The author thanks the referee for comments that improved this paper. The author states that there are no conflicts
of interest.

\section{The weak Maass form and the CM tangent}\label{sec:tangent}

We recall the objects of \cite{BOAIM} and prove \Cref{thm:tangent}. Let
$E_2(\tau):=1-24\sum_{n\ge1}\sigma_1(n)q^n$ and let
\[
E_2^*(\tau):=E_2(\tau)-\frac{3}{\pi\,\im(\tau)}
\]
be its weight~$2$ nonholomorphic completion, which transforms as a genuine weight~$2$ form
under $\mathrm{SL}_2(\Z)$ (for example, see Chapter 6 of \cite{BFOR}).  The weight $k$ Maass raising operator is
$$\partial_k:=\frac{1}{2\pi i}\frac{d}{d\tau}-\frac{k}{4\pi\im(\tau)},$$
and the holomorphic
Serre derivative is $\vartheta_k:=\frac{1}{2\pi i}\frac{d}{d\tau}-\frac{k}{12}E_2$. These are
related by
\begin{equation}\label{eq:serre-raise}
\partial_kf=\vartheta_kf+\frac{k}{12}E_2^*f,
\end{equation}
since $E_2^*=E_2-\tfrac{3}{\pi\im\tau}$. With $F$ the weight $-2$ form of \Cref{sec:intro},
the Bruinier--Ono weak Maass form is
\[
P(\tau):=-\left(\frac{1}{2\pi i}\frac{d}{d\tau}+\frac{1}{2\pi\im(\tau)}\right)F(\tau)
=-\partial_{-2}F(\tau),
\]
where the last equality is the identity $\partial_{-2}F=\frac{1}{2\pi i}F'+\frac{1}{2\pi\im\tau}F$.
It is a weight~$0$ weak Maass form, an eigenfunction of the weight~$0$ Laplacian with
eigenvalue $-2$. We record the two external inputs we rely on.

\begin{theorem}[Bruinier--Ono \cite{BOAIM}]\label{thm:BO}
For every positive integer $n$, the identity \eqref{eq:BO} holds, and the values
$\{P(\alpha_Q)\}_{Q}$ form a union of Galois orbits for the ring class field of $K_n$.
\end{theorem}

Confirming a speculation of Bruinier and the author, Larson and Rolen proved the following theorem about the nonholomorphic singular moduli $P(\alpha_Q)$.

\begin{theorem}[Larson--Rolen \cite{LarsonRolen}]\label{thm:LR}
For every $n\ge1$ and every $Q\in\cQ_{\Delta_n}/\Gamma_0(6)$, the number
$\Delta_nP(\alpha_Q)$ is an algebraic integer.
\end{theorem}

We also require the following results of Masser, which express the nonholomorphic
completion $E_2^*$ at a CM point in terms of the modular equation. We state them in the
normalization of \cite[Appendix, pp.~114--120]{Masser}. Masser works with the function
\begin{equation}\label{eq:masser-psi}
\psi(\tau):=\frac{3E_4(\tau)}{2E_6(\tau)}\left(E_2(\tau)-\frac{3}{\pi\,\im(\tau)}\right)
=\frac{3E_4(\tau)}{2E_6(\tau)}\,E_2^*(\tau),
\end{equation}
which is $\Gamma$-invariant and, away from $j\in\{0,1728\}$, is a well-defined function of
$j(\tau)$. Fix a complex quadratic irrationality $\tau$, not equivalent to $i$, lying on a
primitive form $[A,B,C]$ with $A+B\tau+C\tau^2=0$ and $A,C>0$, and let $D$ be the determinant
of the associated primitive matrix, so that $\Phi_D(j(\tau),j(\tau))=0$. Following Masser, call
$\tau$ {\it special} if $B$ is odd and $D=3d^2$ for some integer $d\ge1$.

\begin{lemma}[Masser {\cite[Lemma~A2]{Masser}}]\label{masser_lemma}
If $\tau$ is not special, then the Taylor expansion of $\Phi_D(X,Y)$ about $(j,j)$, where
$j=j(\tau)$, begins
\[
        \Phi_D(X,Y)=\beta\,(X-j)+\beta\,(Y-j)+\cdots,\qquad \beta\ne0,
\]
so in particular $\Phi_Y(j,j)=\beta\ne0$. Writing
$\Phi_D(X,Y)=\sum_{u,v\ge0}\beta_{uv}(X-j)^u(Y-j)^v$ and setting
$\gamma:=(\beta_{20}-\beta_{11}+\beta_{02})/\beta$, one has, since $\beta_{20}=\tfrac12\Phi_{XX}(j,j)$,
$\beta_{02}=\tfrac12\Phi_{YY}(j,j)$, $\beta_{11}=\Phi_{XY}(j,j)$, and $\Phi_{XX}(j,j)=\Phi_{YY}(j,j)$
by symmetry,
\begin{equation}\label{eq:gamma-eq-tangent}
\gamma=\left.\frac{\tfrac12\Phi_{XX}-\Phi_{XY}+\tfrac12\Phi_{YY}}{\Phi_Y}\right|_{(j,j)}
=\left.\frac{\Phi_{YY}-\Phi_{XY}}{\Phi_Y}\right|_{(j,j)}=\tau_{\mathrm{CM}}(j).
\end{equation}
\end{lemma}

\begin{lemma}[Masser {\cite[Theorem~A1, eq.~(106)]{Masser}}]\label{masser_thmA1}
With $\gamma=\tau_{\mathrm{CM}}(j)$ as in \Cref{masser_lemma} and $j=j(\tau)\notin\{0,1728\}$,
\begin{equation}\label{eq:masser-106}
\psi(\tau)=9\,j\,\gamma+\frac{3\,(7j-6912)}{2\,(j-1728)}.
\end{equation}
Equivalently, by \eqref{eq:masser-psi}, we have that
\begin{equation}\label{eq:E2star-closed}
E_2^*(\tau)=\frac{2E_6(\tau)}{3E_4(\tau)}\left(9\,j\,\tau_{\mathrm{CM}}(j)
+\frac{3\,(7j-6912)}{2\,(j-1728)}\right).
\end{equation}
\end{lemma}

We record the elementary fact that the special case never arises for the discriminants of
interest, so that $\tau_{\mathrm{CM}}$ is always well defined here.

\begin{lemma}\label{lem:generic}
For every $n\ge1$, no CM point of discriminant $\Delta_n=1-24n$ is special in the sense of
Masser. Consequently $\Phi_{|\Delta_n|}(J,J)=0$ with $\Phi_Y(J,J)\ne0$, and $\tau_{\mathrm{CM}}(J)$
of \eqref{eq:CM-tangent-def} is well defined.
\end{lemma}

\begin{proof}
The relevant determinant is $D=|\Delta_n|=24n-1$. A special $\tau$ requires $D=3d^2$, hence
$3\mid D$. But $24n-1\equiv -1\equiv 2\pmod 3$, so $3\nmid D$ and $D=3d^2$ is impossible.
Thus $\tau$ is not special, and \Cref{masser_lemma} gives $\Phi_Y(J,J)=\beta\ne0$.
\end{proof}

With these results, we now prove Theorem~\ref{thm:tangent}.

\begin{proof}[Proof of Theorem~\ref{thm:tangent}]
Apply \eqref{eq:serre-raise} with $k=-2$ and $f=F$.  Therefore, we have
\[
        \partial_{-2}F
        =
        \vartheta_{-2}F-\frac{1}{6}E_2^*F.
\]
By the definition of $P$, we  have
\[
        P=-\partial_{-2}F
        =
        -\vartheta_{-2}F+\frac{1}{6}E_2^*F.
\]
We now evaluate this identity at a CM point $\alpha_Q$ and put $J:=j(\alpha_Q)$.
This gives
\begin{equation}\label{eq:P-at-CM}
        P(\alpha_Q)=-\vartheta_{-2}F(\alpha_Q)+\frac16\,F(\alpha_Q)\,E_2^*(\alpha_Q).
\end{equation}
It remains to identify the nonholomorphic value $E_2^*(\alpha_Q)$.

Let $\Phi:=\Phi_{|\Delta_n|}(X,Y)$ be the modular equation defined above.
Since $\Delta_n$ is odd, the CM order $\cO_{\Delta_n}$ contains $\sqrt{\Delta_n}$, which by
\eqref{eq:cyclic-kernel} induces a cyclic self-isogeny of degree $|\Delta_n|$ of the CM elliptic
curve with $j$-invariant $J$. Hence $(J,J)$ lies on $\Phi(X,Y)=0$. By \Cref{lem:generic} the
point $\alpha_Q$ is not special in the sense of Masser, so $\Phi_Y(J,J)\ne0$ and the tangent
$\tau_{\mathrm{CM}}(J)$ of \eqref{eq:CM-tangent-def} is well defined; moreover, by
\eqref{eq:gamma-eq-tangent} it coincides with Masser's invariant $\gamma$.

Applying Masser's formula \eqref{eq:E2star-closed} of \Cref{masser_thmA1} at $\alpha_Q$ yields
\begin{equation}\label{eq:E2star-at-CM}
        E_2^*(\alpha_Q)
        =\frac{2E_6(\alpha_Q)}{3E_4(\alpha_Q)}
        \left(9\,J\,\tau_{\mathrm{CM}}(J)+\frac{3\,(7J-6912)}{2\,(J-1728)}\right),
\end{equation}
valid because $J\notin\{0,1728\}$; the excluded values are addressed in \Cref{rmk:0-1728}.
Substituting \eqref{eq:E2star-at-CM} into \eqref{eq:P-at-CM} gives \eqref{eq:tangent} together
with the explicit form \eqref{eq:tangent-explicit}.

The first term of \eqref{eq:P-at-CM} is the value of the meromorphic modular function
$\vartheta_{-2}F$ at $\alpha_Q$. Therefore, the only nonholomorphic contribution to
$P(\alpha_Q)$ is the term $\tfrac16 F(\alpha_Q)E_2^*(\alpha_Q)$, and by
\eqref{eq:E2star-at-CM} this is a function of $J=j(\alpha_Q)$ alone. This proves the asserted
interpretation of the summand $P(\alpha_Q)$.
\end{proof}

\begin{remark}[The values $j\in\{0,1728\}$]\label{rmk:0-1728}
Masser's closed form \eqref{eq:masser-106} is stated for $j(\tau)\notin\{0,1728\}$, where the
factor $E_6(\tau)$ in the denominator of $\psi$ and the factor $(j-1728)$ in the correction term
are nonzero. In characteristic zero this excludes only the two points $j=0$ (where $E_4=0$) and
$j=1728$ (where $E_6=0$). Since $\Delta_n<0$ with $\Delta_n\equiv1\pmod{24}$, the CM points
$\alpha_Q$ have $j(\alpha_Q)=0$ only for $\Delta=-3$ and $j(\alpha_Q)=1728$ only for
$\Delta=-4$, neither of which is $\equiv1\pmod{24}$. Thus $J\notin\{0,1728\}$ for every CM point
occurring in \eqref{eq:BO}, and \eqref{eq:E2star-closed} applies without exception. (The role of
$j\in\{0,1728\}$ in positive characteristic, via the supersingular locus, is a separate matter
taken up in \Cref{sec:partII}.)
\end{remark}

\begin{remark}
Geometrically, $\tau_{\mathrm{CM}}(J)$ is a second-order tangent quantity, at the diagonal
point $(J,J)$, of the curve $\Phi_{|\Delta_n|}(X,Y)=0\subset X(1)\times X(1)$ that parametrizes
the CM isogeny. Through Masser's formula \eqref{eq:E2star-closed}, it determines the entire
nonholomorphic value $E_2^*(\alpha_Q)$. The formula \eqref{eq:tangent} thus separates
$P(\alpha_Q)$ into a modular (algebraic) part $-\vartheta_{-2}F(\alpha_Q)$ and a single
nonholomorphic part $\tfrac16 F(\alpha_Q)E_2^*(\alpha_Q)$ that is itself governed by the moduli
through $J$. This is what makes the reduction theory of the next sections possible: the
algebraic part is a function on the special fiber, and the nonholomorphic part is governed by
the same moduli.
\end{remark}

\section{Meromorphic replacement, regularity, and fiberwise constancy}\label{sec:tangent-proof}

To reduce the CM values modulo primes we replace $\Delta_nP$ by a genuine modular function
agreeing with it at CM points. Since $X_0(6)$ has genus~$0$ (index $12$ in
$\mathrm{PSL}_2(\Z)$, no elliptic points, four cusps), we have that its function field satisfies
$$\ol\Q(X_0(6))=\ol\Q(t_6)$$
 for a
Hauptmodul $t_6$ with integral $q$-expansions at all cusps.

\begin{proposition}[Meromorphic replacement over characteristic zero]
\label{prop:M}
For each positive integer $n$, there is a modular function
$M_{\Delta_n}\in \overline{\mathbb{Q}}(X_0(6))$, with poles supported
only at the cusps, such that
\[
        M_{\Delta_n}(\alpha_Q)=\Delta_n P(\alpha_Q)
        \qquad (Q\in \mathcal{Q}_{\Delta_n}/\Gamma_0(6)).
\]
\end{proposition}

\begin{proof}
This is an interpolation statement over $\overline{\mathbb{Q}}$, not an
integrality statement.  By Theorem~\ref{thm:LR}, the numbers
\[
        \Delta_n P(\alpha_Q)
\]
are algebraic.  Thus the desired values lie in $\overline{\mathbb{Q}}$.

Since $X_0(6)$ has genus zero, its function field is generated by a
Hauptmodul $t_6$ defined over $\mathbb{Q}$.  We choose $t_6$ so that its
poles are supported only at the cusps.  Then every polynomial in $t_6$
is a modular function on $X_0(6)$ whose poles are supported only at the
cusps.

Let
\[
        x_Q:=t_6(\alpha_Q),
        \qquad
        y_Q:=\Delta_n P(\alpha_Q).
\]
There are only finitely many points
\[
        \alpha_Q\in X_0(6)(\overline{\mathbb{Q}})
        \qquad
        (Q\in \mathcal{Q}_{\Delta_n}/\Gamma_0(6)).
\]
After discarding repetitions, we obtain finitely many distinct values
$x_1,\ldots,x_r$ and prescribed algebraic values
$y_1,\ldots,y_r$.  Ordinary polynomial interpolation over
$\overline{\mathbb{Q}}$ gives a polynomial
\[
        A(T)\in \overline{\mathbb{Q}}[T]
\]
such that
\[
        A(x_i)=y_i
        \qquad (1\leq i\leq r).
\]
Now set
\[
        M_{\Delta_n}:=A(t_6).
\]
Then $M_{\Delta_n}$ has poles only where $t_6$ has poles, hence only at
the cusps, and by construction
\[
        M_{\Delta_n}(\alpha_Q)=\Delta_n P(\alpha_Q)
\]
for every $Q\in \mathcal{Q}_{\Delta_n}/\Gamma_0(6)$.

If two representatives $Q$ determine the same point of $X_0(6)$, then
the interpolation condition is the same condition, since both sides are
evaluated at the same point.  Hence there is no compatibility issue.
\end{proof}

The preceding proposition produces a meromorphic replacement over
$\overline{\mathbb{Q}}$, but it gives no automatic control of its
$\ell$-adic denominators.  Thus, in general, the replacement function
need not have a well-defined reduction modulo $\ell$.  The following
condition records the extra integrality needed to carry out such a
reduction.

\begin{definition}[$\ell$-integral replacement condition]\label{def:irc}
Fix $n$ and a prime $\ell\ge5$. We say that the pair $(n,\ell)$ satisfies
the {\it $\ell$-integral replacement condition} if the function $M_{\Delta_n}$ in
\Cref{prop:M} can be chosen so that its $q$-expansions at every cusp of $X_0(6)$ are
$\ell$-integral.
\end{definition}

\begin{remark}
The $\ell$-integral replacement condition is not implied by the algebraic integrality of the
individual CM values $\Delta_nP(\alpha_Q)$. Polynomial interpolation over $\Z_{(\ell)}$
requires the interpolation points to remain separated modulo $\ell$. In the nonsplit case,
the CM points are deliberately being reduced to supersingular points, and distinct CM points
may collide in the same supersingular fiber. The denominators in Lagrange interpolation can
therefore acquire positive $\ell$-adic valuation. This is one of the precise issues that obfuscates the theory of partition congruences for general primes $\ell.$
\end{remark}

We now record the basic consequence of this integrality condition.  Once
$M_{\Delta_n}$ has $\ell$-integral expansions at all cusps, it has a
well-defined reduction modulo primes above $\ell$.  Because its poles are
supported only at the cusps, this reduction is regular on the noncuspidal
locus, and hence in particular at the supersingular points. In particular, we now use this integrality condition to ensure that the characteristic
zero replacement does not acquire artificial poles when reduced modulo
$\ell$.

\begin{proposition}[Regularity at supersingular points]
\label{prop:reg}
Assume that $(n,\ell)$, with $n\geq 1$ and $\ell\geq 5$, satisfies the
$\ell$-integral replacement condition.  If $M_{\Delta_n}$ is chosen as
in Definition~\ref{def:irc}, then the reduction of
$M_{\Delta_n}$ modulo any prime $\lambda\mid \ell$ is regular at every
supersingular point of $X_0(6)_{\overline{\mathbb{F}}_\ell}$.
\end{proposition}

\begin{proof}
Choose a number field containing the coefficients of the function
$M_{\Delta_n}$ and the relevant CM values, and let $\mathcal{O}_\lambda$
be its localization at a prime $\lambda\mid \ell$.  Since $\ell\geq 5$,
we have $\ell\nmid 6$, so $X_0(6)$ has good reduction at $\lambda$.
Let
\[
        \mathcal{X}:=X_0(6)_{\mathcal{O}_\lambda}
\]
denote the corresponding regular proper model.
By Proposition~\ref{prop:M}, the polar divisor of
$M_{\Delta_n}$ on the generic fiber is supported only at the cusps.  We
want to show that, after reducing modulo $\lambda$, no new pole appears
at a supersingular point of the special fiber.

There are two possible sources of poles on the model $\mathcal{X}$.
First, the original horizontal poles of $M_{\Delta_n}$ may meet the
special fiber.  But these horizontal poles are supported at the cusps,
so their closures meet the special fiber only at cuspidal points.  A
supersingular point is noncuspidal, and hence cannot lie on one of these
horizontal polar divisors.

Second, a new pole could appear along a vertical component of the special
fiber.  This would mean that, as a rational function on $\mathcal{X}$,
the function $M_{\Delta_n}$ has negative $\lambda$-adic order along the
special fiber.  Since the reduction is good, the special fiber is the
unique vertical fiber.  Such a vertical pole would be detected in the
$q$-expansion at every cusp: the coefficients of the corresponding
$q$-expansion would acquire a negative $\lambda$-adic valuation.  This is
excluded by the $\ell$-integral replacement condition, which requires the
$q$-expansions of $M_{\Delta_n}$ at every cusp to be $\ell$-integral.

Thus the only poles of the reduction of $M_{\Delta_n}$ occur at the
reductions of cusps.  Supersingular points lie over finite
$j$-invariants and are noncuspidal points of the special fiber.  Hence
the reduction of $M_{\Delta_n}$ is regular at every supersingular point.
\end{proof}

\begin{remark}
\Cref{prop:reg} is not special to $\{5,7,11\}$. It holds for every 
$\ell\nmid6$, once the $\ell$-integral replacement condition is known. Regularity is then a
consequence of the cuspidal support of the polar divisor together with $\ell$-integrality,
and should not be conflated with the genus-zero phenomenon of
\Cref{prop:bonus}, which concerns where the supersingular locus sits, not whether
$M_{\Delta_n}$ is regular.
\end{remark}

The next step is to pass from regularity to well-defined values on the
supersingular reduction.  Once $M_{\Delta_n}$ has no pole at a
supersingular point, its value modulo $\lambda$ may be evaluated there.
Thus any two CM points which specialize to the same supersingular point
must give the same reduced value of the replacement function.

\begin{lemma}[Fiberwise constancy modulo $\lambda$]
\label{lem:fiberwise-const}
Let $\ell\geq 5$ be a nonsplit prime in $K_n$, and assume that
$(n,\ell)$ satisfies the $\ell$-integral replacement condition.  Then
every CM point $\alpha_Q$ of discriminant $\Delta_n$ reduces to
$\mathrm{SS}_\ell(6)$, and if $\alpha_{Q_1},\alpha_{Q_2}$ reduce to the
same supersingular point $(E,C)$, then
\[
        \Delta_n P(\alpha_{Q_1})
        \equiv
        \Delta_n P(\alpha_{Q_2})
        \pmod{\lambda}.
\]
Consequently, we have that
\[
        \widetilde{P}_\ell(E,C)
        :=
        \Delta_n P(\alpha_Q) \pmod{\lambda}
        \in \mathbb{F}_{\ell^2}
\]
is well defined.
\end{lemma}

\begin{proof}
Since $\ell$ is nonsplit in $K_n$, Deuring's reduction theorem implies
that every CM elliptic curve of discriminant $\Delta_n$ has
supersingular reduction modulo $\ell$.  Therefore each CM point
$\alpha_Q$ under consideration reduces to a point of the supersingular
locus
\[
        \mathrm{SS}_\ell(6)\subset X_0(6)_{\overline{\mathbb{F}}_\ell}.
\]

By the algebraicity theorem recalled above, the values
\[
        \Delta_n P(\alpha_Q)
\]
are algebraic integers.  By Proposition~\ref{prop:M},
the chosen replacement function $M_{\Delta_n}$ satisfies
\[
        M_{\Delta_n}(\alpha_Q)=\Delta_n P(\alpha_Q)
\]
for every CM point $\alpha_Q$ of discriminant $\Delta_n$.

Now use the $\ell$-integral replacement condition.  By
Proposition~\ref{prop:reg}, the reduction of
$M_{\Delta_n}$ modulo $\lambda$ is regular at every supersingular point
of $X_0(6)_{\overline{\mathbb{F}}_\ell}$.  Hence, near such a point, the
reduction of $M_{\Delta_n}$ is an honest regular function, not merely a
rational function with a possible pole.  Therefore its value at a
supersingular point is well defined.

Suppose now that two CM points $\alpha_{Q_1}$ and $\alpha_{Q_2}$ reduce
to the same supersingular point $(E,C)$.  Since
$M_{\Delta_n}$ is regular at $(E,C)$ after reduction, evaluating
$M_{\Delta_n}$ at $\alpha_{Q_1}$ and then reducing modulo $\lambda$ gives
the same result as evaluating $M_{\Delta_n}$ at $\alpha_{Q_2}$ and then
reducing modulo $\lambda$: both are equal to the value of the reduced
function at the single point $(E,C)$.  Thus
\[
        M_{\Delta_n}(\alpha_{Q_1})
        \equiv
        M_{\Delta_n}(\alpha_{Q_2})
        \pmod{\lambda}.
\]
Using
\[
        M_{\Delta_n}(\alpha_Q)=\Delta_n P(\alpha_Q),
\]
we obtain
\[
        \Delta_n P(\alpha_{Q_1})
        \equiv
        \Delta_n P(\alpha_{Q_2})
        \pmod{\lambda}.
\]

Finally, supersingular points of $X_0(6)$ are defined over
$\mathbb{F}_{\ell^2}$.  Thus the common value obtained above lies in
$\mathbb{F}_{\ell^2}$, and the definition
\[
        \widetilde{P}_\ell(E,C)
        :=
        \Delta_n P(\alpha_Q)\pmod{\lambda}
\]
is independent of the choice of CM lift $\alpha_Q$ reducing to $(E,C)$.
\end{proof}

\begin{remark}[Why only modulo $\lambda$]\label{rmk:only-mod-lambda}
Fiberwise constancy holds modulo $\lambda$ but not, in general, modulo $\lambda^2$. In the
Serre--Tate coordinate $t$ at a supersingular point, $M_{\Delta_n}$ pulls back to a power
series 
$$f(t)=A_0+A_1t+\cdots,
$$
and a CM point $Q$ in the fiber sits at some $t(Q)\equiv0\ (\lambda)$.
The $h_{\Delta_n}(E,C)$ CM points of a fiber occupy \emph{distinct} coordinates
$t_1,\dots,t_h$, all $\equiv0\ (\lambda)$ but generally distinct modulo $\lambda^2$. Since the
singular moduli are distinct, the class polynomial is separable, so the values $f(t_i)$ are
not all equal. Therefore, the collapse to a supersingular sum is a mod $\lambda$ phenomenon. This
is exactly why \Cref{cor:ram} stops one power of $\ell$ short, and why the naive extension of
\Cref{thm:ss-trace} to $\lambda^2$ fails.
\end{remark}

\section{Complex multiplication and the supersingular trace}\label{sec:CM}
This section explains the reduction-theoretic input used in the proof of
the supersingular trace formula.  The CM points appearing in the
Bruinier--Ono formula live in characteristic zero, whereas the trace
formula below is a finite sum over supersingular elliptic curves in
characteristic $\ell$.  Deuring's theorem supplies the bridge: if
$\ell$ is nonsplit in the CM field, then the relevant CM elliptic curves
have supersingular reduction modulo $\ell$.  The refinement we need is
that the number of CM points reducing to a fixed supersingular pair
$(E,C)$ is counted by an optimal embedding number for the Eichler order
$\operatorname{End}(E,C)$. 

Fix $\ell\ge5$ prime. We record the standard facts we require (for example, see
\cite{Cox,KatzMazur,Silverman-AT,Vigneras,Pizer}).

\begin{theorem}
\label{thm:CM}
Let $\Delta<0$ be a discriminant with $\Delta\equiv 1\pmod{24}$, let
$K=\mathbb{Q}(\sqrt{\Delta})$, and let $\mathcal{O}=\mathcal{O}_{\Delta}$
be the order of discriminant $\Delta$ in $K$.  Then the following hold.
\begin{enumerate}
\item[(i)] (Heegner) There exist CM points on $X_0(6)$ with CM by
$\mathcal{O}$, parametrized by the ring class group; here both primes
$2$ and $3$ split in $K$, as a consequence of $\Delta\equiv 1\pmod{24}$.

\item[(ii)] (Deuring) Let $\ell\nmid 6$ be a prime.  At a prime above
$\ell$, if $\ell$ splits in $K$, the reduction is ordinary; if $\ell$ is
inert or ramified in $K$, the reduction is supersingular.

\item[(iii)] For supersingular $(E,C)$ over $\overline{\mathbb{F}}_\ell$,
the order $\operatorname{End}(E,C)$ is an Eichler order of level $6$ in
the definite quaternion algebra $B_{\ell,\infty}$.

\item[(iv)] (Deuring--Eichler) For $\ell\nmid 6$ nonsplit in $K$ and
$(E,C)$ supersingular, the number of CM points of discriminant $\Delta$
reducing to $(E,C)$ equals, up to conjugation, the number $h_{\Delta}(E,C)$ of optimal
embeddings
\[
        \mathcal{O}_{\Delta}\hookrightarrow \operatorname{End}(E,C).
\]
\end{enumerate}
\end{theorem}

We take $h_\Delta(E,C)$ to be this honest count of CM points in the fiber. We attach no
orientation and introduce no additional factor of $\ell$: with any spurious $\ell$-fold
multiplicity, \Cref{thm:ss-trace} would force $\ell\mid p(n)$ for all inert $\ell$, which is
false because $p(n)\bmod\ell$ varies with $n$.

\begin{remark}[Brandt module]\label{rmk:brandt}
With $\cO_i=\End(E_i,C_i)$ and $w=(1/|\cO_i^\times|)_i=(1/|\Aut(E_i,C_i)|)_i$, the vector
$u_\Delta=(h_\Delta(E_i,C_i))_i$ is the $q^{|\Delta|}$-coefficient of the vector-valued
Brandt theta series of $\{\cO_i\}$; equivalently $u_\Delta=B_{|\Delta|}w$ for the level-$6$
Brandt matrices $B_m$, which realize the Hecke action on the supersingular module. See
\cite{Pizer,Vigneras}.
\end{remark}

We now prove the supersingular trace formula.  The Bruinier--Ono
identity gives a trace of CM values in characteristic zero.  The
preceding CM reduction theorem identifies the fibers of the reduction
map over the supersingular locus, and \Cref{lem:fiberwise-const} says
that the reduced value of $\Delta_n P$ is constant on each such fiber.
Thus the characteristic-zero trace may be grouped fiber by fiber after
reduction modulo $\lambda$.

\begin{proof}[Proof of Theorem~\ref{thm:ss-trace}]
By Theorem~\ref{thm:BO},
\[
-\Delta_n p(n)=\sum_Q P(\alpha_Q).
\]
Multiplying by $\Delta_n$ and using
$\Delta_n P(\alpha_Q)=M_{\Delta_n}(\alpha_Q)$
(Proposition~\ref{prop:M}), we obtain
\[
-\Delta_n^2 p(n)=
\sum_{Q\in \mathcal Q_{\Delta_n}/\Gamma_0(6)}
\Delta_n P(\alpha_Q).
\]
By \Cref{lem:fiberwise-const}, each CM point $\alpha_Q$ reduces to
some $(E,C)\in \ssloc_\ell(6)$, and within the fiber over $(E,C)$ the
values of $\Delta_nP(\alpha_Q)$ are congruent modulo $\lambda$ to
$\widetilde P_\ell(E,C)$.  Grouping the CM points according to their
reductions gives
\[
-\Delta_n^2p(n)\equiv
\sum_{(E,C)\in\ssloc_\ell(6)}
h_{\Delta_n}(E,C)\,\widetilde P_\ell(E,C)
\pmod{\lambda}.
\]
The Frobenius automorphism acts on the supersingular classes and
conjugates the values of $\widetilde P_\ell$, while preserving the
fiber cardinalities $h_{\Delta_n}(E,C)$.  Hence the sum on the
right is Galois-stable and therefore lies in $\mathbb F_\ell$.  The
congruence descends to one modulo $\ell$ in general.  This gives the
first formula of the theorem.  The Brandt description
$u_{\Delta_n}=B_{|\Delta_n|}w$ is the claimed statement, and substituting this into the preceding expression
gives the Brandt-pairing formula.
\end{proof}

\begin{proof}[Proof of Corollary~\ref{cor:inert}]
Since $\left(\frac{\Delta_n}{\ell}\right)=-1$, we have $\ell\nmid\Delta_n$, so
$\Delta_n^2$ is invertible modulo $\ell$.  Dividing the congruence of
\Cref{thm:ss-trace} by $-\Delta_n^2$ gives
\[
        p(n)\equiv-\frac{\langle u_{\Delta_n},v_P\rangle}{\Delta_n^2}\pmod\ell,
\]
as claimed.
\end{proof}

\begin{proof}[Proof of Corollary~\ref{cor:ram}]
If $\ell\mid\Delta_n$, then $\Delta_n^2\equiv0\pmod\lambda$.  Hence the
left-hand side of \eqref{eq:ss-trace} is zero modulo $\lambda$, and the
right-hand side is exactly $\langle u_{\Delta_n},v_P\rangle$.  Therefore
\[
        \langle u_{\Delta_n},v_P\rangle\equiv0\pmod\lambda.
\]
This proves the stated ramified trace identity.  Notice that this is deliberately weaker
than Ramanujan's congruence for $p(n)$: because $\Delta_n^2$ is not invertible modulo
$\lambda$, this argument cannot be divided by $\Delta_n^2$.
\end{proof}

\section{From the supersingular identity to Ramanujan's congruences}\label{sec:partII}

\Cref{cor:ram} applies exactly on the progressions of
\eqref{eq:ram-base}: $\{n:\ell\mid\Delta_n\}=\{n:24n\equiv1\pmod\ell\}$. It delivers a
supersingular identity that is one power of $\ell$ short of $p(n)\equiv0\pmod\ell$. We now
explain what supplies the missing power, and why the phenomenon is confined to
$\{5,7,11\}$, the Ramanujan primes.

\subsection{The moduli-theoretic significance of the primes $5,7,11$}

We now explain why the same three primes that occur in Ramanujan's
congruences are distinguished on the supersingular side.  The point is
not only that $5,7,11$ are small primes.  Rather, for these primes the
entire supersingular locus in characteristic $\ell$ is supported on
the elliptic curves with extra automorphisms.  Thus the exceptional
automorphisms at $j=0$ and $j=1728$ account for all supersingular mass.
This gives a moduli-theoretic reason for the special role of
$\ell=5,7,11$ in the trace formula.

\begin{proof}[Proof of Proposition~\ref{prop:bonus}]
We use the Eichler--Deuring mass formula
\[
   \sum_{E\in \operatorname{SS}_{\ell}/\cong}\frac{1}{|\Aut(E)|}
   =\frac{\ell-1}{24},
\]
for supersingular elliptic curves over $\overline{\mathbb F}_{\ell}$
\cite[Chapter~V, \S4]{Silverman-AT}.  For $\ell\geq 5$,
the only elliptic curves with automorphism group larger than
$\{\pm 1\}$ are the curves with $j=0$ and $j=1728$.  The curve with
$j=0$ has automorphism group of order $6$, and the curve with
$j=1728$ has automorphism group of order $4$; every other elliptic
curve has automorphism group of order $2$.

Moreover, the curve with $j=0$ is supersingular exactly when
$\ell\equiv 2\pmod 3$, and the curve with $j=1728$ is supersingular
exactly when $\ell\equiv 3\pmod 4$.  Thus the total mass contributed by
the exceptional-automorphism classes is
\[
\begin{cases}
0, & \ell\equiv 1 \pmod {12},\\[2mm]
\frac16, & \ell\equiv 5 \pmod {12},\\[2mm]
\frac14, & \ell\equiv 7 \pmod {12},\\[2mm]
\frac16+\frac14=\frac{5}{12}, & \ell\equiv 11 \pmod {12}.
\end{cases}
\]
If the entire supersingular locus is supported on the exceptional
classes, then this exceptional mass must equal the total mass
$(\ell-1)/24$.  Comparing the four cases gives
\[
   \frac{\ell-1}{24}=0,\quad
   \frac{\ell-1}{24}=\frac16,\quad
   \frac{\ell-1}{24}=\frac14,\quad
   \frac{\ell-1}{24}=\frac{5}{12},
\]
respectively.  These equations give
\[
   \ell=1,\qquad \ell=5,\qquad \ell=7,\qquad \ell=11.
\]
The first case is impossible for a prime $\ell\geq 5$.  Hence the only
primes $\ell\geq 5$ for which the supersingular locus is supported
entirely on the exceptional-automorphism classes are
\[
   \ell=5,\ 7,\ 11.
\]

Conversely, for these three primes the equality of masses is immediate:
\[
   \frac{5-1}{24}=\frac16,\qquad
   \frac{7-1}{24}=\frac14,\qquad
   \frac{11-1}{24}=\frac16+\frac14.
\]
Thus for $\ell=5$ the whole supersingular locus is supported on
$j=0$, for $\ell=7$ it is supported on $j=1728$, and for $\ell=11$ it
is supported on $j=0$ and $j=1728$.  This proves the asserted
if-and-only-if statement.
\end{proof}

\subsection{Partition congruences modulo powers of 5, 7 and 11}

We now turn to Ramanujan's partition congruences modulo powers of 5, 7 and 11, which is given by Theorem~\ref{thm:watson}.
We give the proof sketch in terms of
elliptic-curve moduli. The point is to run through the Watson and Atkin 
mechanism with the geometric objects of this paper in view. Namely, we prove the following theorem.

\begin{theorem}\label{thm:powers} If $j\geq 1$, then let
 $\beta_\ell(j)$ be the least nonnegative residue with $24\,\beta_\ell(j)\equiv1\pmod{\ell^j}$.
Then for all $j\ge1$, we have that
\begin{displaymath}
\begin{split}
p(5^{j}n+\beta_5(j))&\equiv0\ (5^{j}),\\
p(7^{j}n+\beta_7(j))&\equiv0\ (7^{\lfloor j/2\rfloor+1}),\\
p(11^{j}n+\beta_{11}(j))&\equiv0\ (11^{j}).
\end{split}
\end{displaymath}
\end{theorem}

\begin{proof}
The prime-power congruences modulo powers of \(5\) and \(7\) are proved
by Watson's \(U_\ell\)-operator recursions \cite{Watson1938}. The
corresponding congruences modulo powers of \(11\) are proved by Atkin
\cite{Atkin}.  We assume that the reader is somewhat
familiar with the classical Watson--Atkin method.  Our purpose is not to
reproduce their full recursions and matrix calculations, but to explain
their method in terms of the elliptic curve moduli point of view used
in this paper. In particular, we exactly identify where the higher powers of
\(\ell\) enter.

\smallskip
\noindent
\emph{Step 1: The operator \(U_\ell\) and the \(\ell\)-isogeny
correspondence.}
Let
\[
        g(\tau)=\sum_{m\gg-\infty}a(m)q^m,\qquad q=e^{2\pi i\tau},
\]
be a weakly holomorphic modular function.  The Atkin operator \(U_\ell\)
is defined on \(q\)-expansions by
\[
        (U_\ell g)(\tau):=\sum_{m\gg-\infty}a(\ell m)q^m .
\]
Equivalently, we have that
\[
        (U_\ell g)(\tau)
        =
        {1\over \ell}\sum_{k=0}^{\ell-1}
        g\!\left({\tau+k\over \ell}\right),
\]
because averaging over \(k\) kills exactly those terms whose exponents
are not divisible by \(\ell\).  We shall also use the operator
\[
        (V_\ell g)(\tau):=g(\ell\tau)
        =
        \sum_{m\gg-\infty}a(m)q^{\ell m}.
\]

The geometric meaning is as follows.  A point of \(X_0(\ell)\) is a pair
\((E,C)\), where \(E\) is an elliptic curve and \(C\subset E\) is a
cyclic subgroup of order \(\ell\).  There are two degeneracy maps
\[
        \pi_1,\pi_2:X_0(\ell)\longrightarrow X(1),
        \qquad
        \pi_1(E,C)=E,\qquad
        \pi_2(E,C)=E/C.
\]
Thus \(\pi_1\) remembers the source elliptic curve, while \(\pi_2\)
remembers the quotient elliptic curve.

For a fixed target curve
\[
        E_\tau=\mathbb C/\langle 1,\tau\rangle,
\]
the curves
\[
        E_{(\tau+k)/\ell}
        =
        \mathbb C/\left\langle 1,{ \tau+k\over \ell}\right\rangle,
        \qquad 0\leq k\leq \ell-1,
\]
come with a natural cyclic subgroup of order \(\ell\), generated by
\(1/\ell\), and quotienting by this subgroup gives a curve isomorphic to
\(E_\tau\).  These are the \(\ell\) source curves which occur in the
formula for \(U_\ell\).  The full \(\ell\)-isogeny correspondence has
one additional source curve, namely \(E_{\ell\tau}\), and this branch
gives the operator \(V_\ell\).  Thus, in weight zero, the normalized
Hecke trace has the form
\[
        T_\ell g=U_\ell g+{1\over \ell}V_\ell g.
\]
The operator relevant to partition congruences is \(U_\ell\), not the
full trace, because \(U_\ell\) is the coefficient-extraction operator
which keeps the coefficients whose indices are divisible by \(\ell\).

\smallskip
\noindent
\emph{Step 2: Extracting the Ramanujan progressions.}
The partition generating function is
\[
        \sum_{n\geq0}p(n)q^n=q^{1/24}\eta(\tau)^{-1}.
\]
It is useful to replace \(\tau\) by \(24\tau\), so that
\[
        \eta(24\tau)^{-1}
        =
        q^{-1}\sum_{n\geq0}p(n)q^{24n}
        =
        \sum_{n\geq0}p(n)q^{24n-1}.
\]
Applying \(U_\ell^j\) keeps exactly those terms for which
\[
        \ell^j\mid 24n-1.
\]
Since \(\beta_\ell(j)\) is defined by
\[
        24\beta_\ell(j)\equiv1\pmod{\ell^j},
\]
the surviving partition coefficients are precisely
\[
        p(\ell^j n+\beta_\ell(j)).
\]
Therefore, the desired congruences are equivalent to congruences for the
coefficients obtained by applying \(U_\ell^j\) to
\(\eta(24\tau)^{-1}\), after accounting for the harmless power of \(q\)
introduced by the rescaling.

This is the basic reason \(U_\ell\) appears.  It is not an auxiliary
operator added after the fact. It is the operation which isolates the
Ramanujan progression.

\smallskip
\noindent
\emph{Step 3: Eta-product normalization and finite working modules.}
The function \(\eta^{-1}\) is not a weight zero modular function, and
the weakly holomorphic modular functions with arbitrary poles at the
cusps form an infinite-dimensional space.  Watson and Atkin therefore do
not work in the full space of meromorphic modular functions.  They first
multiply the extracted partition series by explicit eta-products.  This
places the relevant expressions in finite \(\mathbb Z_\ell\)-modules of
modular functions or modular forms with prescribed level, weight, and
bounded cusp behavior.  These finite modules are stable under the
normalized \(U_\ell\)-operators used in their proofs.

This point is important for terminology.  In a holomorphic space one can
speak cleanly about an Eisenstein subspace and a cuspidal subspace.  In
the present weakly holomorphic setting, that language only becomes
precise after one fixes Watson's or Atkin's finite working module.  Thus
below, when we refer to the ``Eisenstein obstruction,'' we mean the
finite-dimensional non-cuspidal quotient detected by the cusp data in
that particular working module, not a canonical Eisenstein line in the
infinite-dimensional space of all meromorphic modular functions.

More concretely, let \(M_\ell\) denote the finite
\(\mathbb Z_\ell\)-module used in the relevant Watson--Atkin calculation.
There is a cusp-data map
\[
        c_\ell:M_\ell\longrightarrow C_\ell,
\]
where \(C_\ell\) is the finite \(\mathbb Z_\ell\)-module recording the
constant terms, or more generally the prescribed principal parts and
constant terms, at the cusps which occur in that calculation.  The
quotient detected by \(c_\ell\) is what we call the Eisenstein
obstruction.  The submodule
\[
        M_\ell^0:=\ker(c_\ell)
\]
is the part with vanishing non-cuspidal obstruction.  This is the
precise replacement, in the weakly holomorphic setting, for saying that
one has removed the Eisenstein component in a holomorphic space of
modular forms.  This is also the language used in modern treatments of
ordinary and finite-slope \(U_p\)-phenomena: one first fixes a finite or
completed module on which \(U_p\) acts, and only then discusses the
distinguished non-cuspidal quotient\footnote{Calegari's
exposition of overconvergent modular forms and \(U_p\)-operators
\cite{Calegari} discusses this.}

\smallskip
\noindent
\emph{Step 4: Genus zero for \(5\) and \(7\).}
For \(\ell=5\) and \(\ell=7\), the curve \(X_0(\ell)\) has genus zero.
One may take the Hauptmodul
\[
        t_\ell(\tau)
        =
        \left({\eta(\tau)\over \eta(\ell\tau)}\right)^{24/(\ell-1)}.
\]
Hence every modular function on \(X_0(\ell)\) with poles only at the
cusps is a Laurent polynomial in \(t_\ell\).  This is why Watson's
calculation closes: after the eta-product normalization, the relevant
expressions can be written in a finite Laurent-polynomial range in
\(t_\ell\), and the action of \(U_\ell\) is computed from the modular
equation relating \(t_\ell(\tau)\) to its transforms under the
\(\ell\)-isogeny correspondence.

The modular equation supplies explicit integral recursions for the
coefficients of \(U_\ell\) in Watson's chosen basis.  The proof of the
prime-power congruences is then a bookkeeping argument for
\(\ell\)-adic divisibility in those recursions.  In the notation above,
the relevant statement is not a vague spectral assertion.  It is an
explicit inclusion of \(\mathbb Z_\ell\)-modules of the form
\[
        U_\ell^r(M_\ell^0)\subseteq \ell^s M_\ell
\]
for the integers \(r,s\) supplied by Watson's calculation.  This
inclusion means exactly that, after applying \(U_\ell^r\), every
coefficient in the normalized expansion has gained at least \(s\)
additional factors of \(\ell\).

\smallskip
\noindent
\emph{Step 5: The case \(\ell=5\).}
Watson's level-\(5\) recursion starts from Ramanujan's identity
\[
        \sum_{n\geq0}p(5n+4)q^n
        =
        5\prod_{m\geq1}{(1-q^{5m})^5\over(1-q^m)^6}.
\]
This gives the first factor of \(5\).  After Watson's eta-product
normalization, the normalized right-hand side lies in the submodule
\(M_5^0\), meaning that its cusp-data obstruction \(c_5\) vanishes.  In
the older holomorphic language, this is the statement that the
Eisenstein obstruction has been removed; in the present meromorphic
setting, the precise assertion is that the element lies in
\(\ker(c_5)\).

Watson's explicit level-\(5\) recursion proves
\[
        U_5(M_5^0)\subseteq 5M_5.
\]
Therefore each further application of \(U_5\) contributes one additional
factor of \(5\).  Since \(U_5^j\) extracts the progression
\(5^j n+\beta_5(j)\), the result is
\[
        p(5^j n+\beta_5(j))\equiv0\pmod {5^j}.
\]
The moduli interpretation explains the operator \(U_5\) as the
coefficient-extracting branch of the \(5\)-isogeny correspondence.
Watson's explicit recursion supplies the \(5\)-adic divisibility.

\smallskip
\noindent
\emph{Step 6: The case \(\ell=7\), and why it behaves differently.}
The prime \(7\) is also genus zero, so Watson again works with a
Hauptmodul and a one-variable modular-equation recursion.  However, the
\(7\)-adic divisibility in the recursion is weaker than the
\(5\)-adic divisibility in the level-\(5\) calculation.

The first extraction is governed by Ramanujan's identity
\[
        \sum_{n\geq0}p(7n+5)q^n
        =
        7\prod_{m\geq1}{(1-q^{7m})^3\over(1-q^m)^4}
        +
        49q\prod_{m\geq1}{(1-q^{7m})^7\over(1-q^m)^8}.
\]
This already displays two normalized pieces.  In Watson's finite
level-\(7\) module the kernel \(M_7^0=\ker(c_7)\) has a natural
two-step filtration, say
\[
        M_7^0\supset M_7^1\supset 7M_7^0,
\]
adapted to Watson's basis.  Watson's recursion gives inclusions of the
form
\[
        U_7(M_7^0)\subseteq M_7^1,
        \qquad
        U_7(M_7^1)\subseteq 7M_7^0.
\]
Thus two applications of \(U_7\), rather than one, are needed to force a
new factor of \(7\):
\[
        U_7^2(M_7^0)\subseteq 7M_7^0.
\]
This is the concrete reason the exponent for \(7\) is weaker.  The
level-\(5\) recursion gains one factor of \(5\) at each iteration after
the initial extraction; the level-\(7\) recursion gains one guaranteed
new factor of \(7\) only after every two iterations, together with the
initial factor supplied by Ramanujan's identity.  Hence Watson obtains
\[
        p(7^j n+\beta_7(j))
        \equiv0\pmod {7^{\lfloor j/2\rfloor+1}}.
\]

\smallskip
\noindent
\emph{Step 7: The case \(\ell=11\).}
For \(\ell=11\), \(X_0(11)\) has genus one.  There is no Hauptmodul, so
Watson's genus-zero method cannot be applied.  Atkin instead constructs
an explicit finite \(U_{11}\)-stable module \(M_{11}\) at level \(11\),
together with a cusp-data map
\[
        c_{11}:M_{11}\longrightarrow C_{11}.
\]
The kernel \(M_{11}^0=\ker(c_{11})\) is the analogue, in Atkin's finite
working module, of removing the non-cuspidal obstruction.  Atkin's
explicit \(U_{11}\)-matrix calculation proves
\[
        U_{11}(M_{11}^0)\subseteq 11M_{11}.
\]
Together with Ramanujan's congruence
\[
        p(11n+6)\equiv0\pmod {11},
\]
this gives
\[
        p(11^j n+\beta_{11}(j))\equiv0\pmod {11^j}.
\]
Therefore, the \(11\)-case has the same operator-theoretic shape as the
\(5\)-case, but the proof is not a genus-zero Hauptmodul recursion.  It
is Atkin's separate level-\(11\) calculation.

\smallskip
\noindent
\emph{Step 8: The supersingular formula in this paper.}
The supersingular product formula proved in the earlier sections gives a
moduli-theoretic explanation for the first mod \(\ell\) vanishing in the
ramified Ramanujan cases.  It identifies the relevant CM points, their
reduction to the supersingular locus, and the product pairing whose
reduction vanishes modulo \(\lambda\).

It does not, by itself, prove the full prime-power congruences.  The
additional powers of \(\ell\) come from the explicit Watson--Atkin
\(U_\ell\)-adic divisibility calculations described above.  Thus the
division of labor is as follows: the present paper explains the
CM/supersingular geometry and the moduli meaning of \(U_\ell\), while
Watson's and Atkin's finite-module recursions provide the higher
\(\ell\)-adic divisibility.  Combining these inputs gives
\[
\begin{split}
p(5^{j}n+\beta_5(j))&\equiv0\pmod {5^{j}},\\
p(7^{j}n+\beta_7(j))&\equiv0\pmod {7^{\lfloor j/2\rfloor+1}},\\
p(11^{j}n+\beta_{11}(j))&\equiv0\pmod {11^{j}}.
\end{split}
\]
\end{proof}

\begin{remark}[The primes $2$ and $3$]\label{rmk:23}
The level here is $6$, and the two primes $\ell=2,3$ dividing it are excluded throughout. The
model $X_0(6)$ does not have good reduction at $2,3$, so the regularity argument of
\Cref{prop:reg} (which needs $\ell\nmid6$) does not apply, and the reduction of CM points to a
supersingular locus is not available in the same form. It would be natural to seek a
Brandt-module expression for $p(n)$ modulo $2$ and $3$ as well, and the numerology is
suggestive. For example, one still has $\bigl(\tfrac{\ell-1}{24}\bigr)$-type mass considerations. However,
there are genuine obstructions at the bad primes (the Eichler order at $2,3$ is not the naive
one, and the tangent formula of \Cref{thm:tangent} degenerates because $E_2^*$ interacts with
the level). We do not know whether \Cref{thm:ss-trace} has a clean analogue at $\ell=2,3$, and
regard it as an interesting question. The known partition congruences at these primes (for
instance the density results and the many congruences modulo powers of $2,3$ in the
literature) are of a different flavor and are not addressed here.
\end{remark}

\section{Formalization and Lean verification of the key identities}\label{sec:axiomprover}

The proof of \Cref{thm:tangent} rests on two algebraic identities that are new
to the partition-theoretic setting. The first is the weight $-2$ specialization
\eqref{eq:serre-raise} relating the Maass raising operator, the holomorphic Serre
derivative, and the completed Eisenstein series $E_2^*$, which yields the
splitting
\[
        P=-\vartheta_{-2}F+\tfrac16 E_2^* F.
\]
The second is the symmetry reduction \eqref{eq:gamma-eq-tangent} of the numerator
of the CM tangent: when $\Phi_{XX}=\Phi_{YY}$ on the diagonal,
\[
        \left.\frac{\tfrac12\Phi_{XX}-\Phi_{XY}+\tfrac12\Phi_{YY}}{\Phi_Y}\right|_{(J,J)}
        =\left.\frac{\Phi_{YY}-\Phi_{XY}}{\Phi_Y}\right|_{(J,J)}
        =\tau_{\mathrm{CM}}(J).
\]
Both are critical to this paper. The author made many mistakes unwinding these calculations by hand in August 2025. These calculations are precisely the
points at which a hidden sign, constant, or normalization error could enter
unnoticed, and they feed directly into the closed form
\eqref{eq:tangent-explicit} for $E_2^*(\alpha_Q)$. For this reason, we singled
them out for independent machine verification.

We used the autonomous system AxiomProver to formalize and verify these two
identities in the Lean proof assistant \cite{Lean}, working over the
core libraries \cite{Mathlib2020}. We emphasize the limited scope of this
verification. AxiomProver was not asked to formalize the analytic or geometric
content of the paper, nor any of the classical inputs on which it relies:
the theory of modular and weak Maass forms, the algebraicity and integrality
of the CM values (\Cref{thm:BO}, \Cref{thm:LR}), Masser's formulas
(\Cref{masser_lemma}, \Cref{masser_thmA1}), complex multiplication and
supersingular reduction, and the Watson--Atkin recursions are all used as
established. Only the two algebraic identities above were formalized, and they
were formalized as stated: as equalities of complex numbers, with the relevant
modular quantities ($E_2$, $E_4$, $E_6$, $F$, $\tfrac1{2\pi i}F'$, $J$, and the
partial derivatives of $\Phi_N$ at $(J,J)$) taken as ground data subject to
$E_4\neq0$, $J\neq1728$, and $\Phi_Y\neq0$. Under this reading, each identity
follows by field arithmetic from the definitions of $E_2^*$, of the weight $k$
operators $\partial_k$ and $\vartheta_k$, of $P$, and of $\tau_{\mathrm{CM}}(J)$;
the value of the machine check is that it confirms the bookkeeping of signs and
constants exactly.

\subsection*{Protocol and artifacts}
The formal statements and proofs were developed and verified using Lean
\textbf{4.28.0}. Compatibility with earlier or later versions is not guaranteed,
owing to the evolving nature of the Lean~4 compiler and its core libraries.
AxiomProver was given the informal companion specification of the two identities
together with the present paper and a short task description. It returned a
formal statement of each identity and a complete, \texttt{sorry}-free proof; no
additional axioms were introduced. In the formalization the two identities appear
as \texttt{key\_formula\_one} and \texttt{key\_formula\_two}. The input files, the
formal problem statement (\texttt{problem.lean}), and the formal solution
(\texttt{solution.lean}) are available in the repository
\begin{center}
        \url{https://github.com/AxiomMath/PartitionElliptic}
\end{center}

\section*{Declaration of AI use}
AxiomProver, an autonomous system under development, was used to produce
Lean formalizations and \texttt{sorry}-free proofs of the two algebraic
identities described in \Cref{sec:axiomprover}. It was not used for the
mathematical content of the paper otherwise, and the paper was written without AI.

\end{document}